\documentclass[journal]{IEEEtran}
\usepackage{algorithm}
\usepackage{algorithmic}
\usepackage{todonotes}
\usepackage{float}
\usepackage{breqn}
\usepackage{caption}
\usepackage{subcaption}
\usepackage{makecell}
\usepackage{multirow}
\usepackage{cite}

\ifCLASSINFOpdf
\else
\fi
\usepackage{amsmath}
\usepackage{graphicx}

\begin{document}
%
\title{Enhancing Time Series Aggregation For Power System Optimization Models: Incorporating Network and Ramping Constraints}
%
%

\author{David~Cardona-Vasquez,
        Thomas~Klatzer,
        Sonja~Wogrin
}

%
%

\markboth{Journal of \LaTeX\ Class Files,~Vol.~14, No.~8, August~2015}%
{Shell \MakeLowercase{\textit{et al.}}: Bare Demo of IEEEtran.cls for IEEE Journals}
%



\maketitle

\begin{abstract}

*Power system optimization models are large mathematical models used by researchers and policymakers that pose tractability issues when representing real-world systems. Several aggregation techniques have been proposed to address these computational challenges and it remains a relevant topic in power systems research. In this paper, we extend a recently developed Basis-Oriented time series aggregation approach used for power system optimization models that aggregates time steps within their Simplex basis. This has proven to be an exact aggregation for simple economic dispatch problems. We extend this methodology to include network and ramping constraints; for the latter (and to handle temporal linking), we develop a heuristic algorithm that finds an exact partition of the input data, which is then aggregated. Our numerical results, for a simple 3-Bus system, indicate that: with network constraints only, we can achieve a computational reduction by a factor of 1747 (measured in the number of variables of the optimization model), and of 12 with ramping constraints. Moreover, our findings indicate that with temporal linking constraints, aggregations of variable length  must be employed to obtain an exact result (the same objective function value in the aggregated model) while maintaining the computational tractability, this implies that the duration of the aggregations does not necessarily correspond to commonly used lengths like days or weeks. Finally, our results support previous research concerning the importance of extreme periods on model results.

\end{abstract}

\begin{IEEEkeywords}
power system optimization, mathematical modeling, dimensionality reduction, renewable energy sources, time series aggregation, linear programming.
\end{IEEEkeywords}

%
\IEEEpeerreviewmaketitle

\subsection{Sets and Subsets}
\begin{IEEEdescription}[\IEEEusemathlabelsep\IEEEsetlabelwidth{$a(m,n,l)$}]
\item[$k$] Set of hours
\item[$r$] Set of representative periods
\item[$l$] Set of lines
\item[$g$] Set of generation units
\item[$w, t \in g$] Generation units (wind and thermal)
\item[$i, j$] Sets of nodes
\end{IEEEdescription}

\subsection{Variables and Parameters}
\begin{IEEEdescription}[\IEEEusemathlabelsep\IEEEsetlabelwidth{$a(m,n,l)$}]
\item[$p_{g,k}$] Power produced by $w$ (or $t$) during hour $k$ in MWh
\item[$nsp_{k}$] Non-supplied power during hour $k$ in MWh
\item[$nsp_{k,i}$] Non-supplied power during hour $k$ in MWh at node $i$
\item[$f_{k,i, j}$] Flow from node $i$ to $j$ during hour $k$ in MWh
\item[$\overline{F_{l}}$] Flow limit for line $l$ in MW
\item[$\overline{P_{g,k}}$] Installed capacity of $g$ in MW
\item[$\underline{P_{g,k}}$] Minimum production for $g$ in MW
\item[$CF_{k}$] Capacity factor at time $k$ 
\item[$D_{k/r}$] Power demand during hour $k$ or representative period $r$ in MW
\item[$C^{g}$] Variable cost of $g$ in \texteuro/MWh
\item[$C^{nsp}$] Cost of $nsp$ in \texteuro/MWh
\item[$C^{N}$] Transmission cost in \texteuro/MWh
\item[$RD_t$] Maximum ramping down of $t$ in MW
\item[$RU_t$] Maximum ramping up of $t$ in MW
\item[$VC_{g}$] Variable cost of $g$ in \texteuro/MWh
\item[$MC_{k}$] Marginal cost at time $k$ in \texteuro/MWh
\end{IEEEdescription}

\section{Intrdouction}
\label{sec:intro}

\IEEEPARstart{P}{ower} System Optimization Models (PSOMs) are widely used for planning and policy-making toward sustainable and clean energy systems. However, due to real power systems' spatio-temporal size and technical complexity, PSOMs can result in computationally intractable problems. Computational intractability in PSOMs stems from the multiple dimensions they try to represent, e.g., the technical details, the uncertainty concerning the system's demand and generators' availability, the spatial representation, or the granularity and length of the time horizon, etc. To overcome this, modelers apply aggregation techniques to approximate full PSOMs with PSOMs of reduced size to derive practical results within reasonable CPU times. One subset of such techniques is time series aggregation (TSA); which aims to replace a full hourly or even sub-hourly PSOM with a smaller model using a simplified time dimension, allowing for faster model runs, and maintaining, at least to some degree, the accuracy of the results. There are many methods to achieve this, as reviewed in \cite{en13030641} and \cite{TEICHGRAEBER2022111984}, and each approach comes with its advantages and shortcomings.

In the literature \cite{en13030641}, TSA techniques are subdivided in a-priori and a-posteriori methods. A-priori methods rely on the input space, e.g., demand time series, capacity factors of renewables etc., as most of the models cannot be run without performing an aggregation procedure; so, a-priori techniques, have little regard for the optimization model performance or its structure, and even today, they remain as state-of-the-art TSA techniques for PSOMs\cite{Sarajpoor2021} \cite{Liu2018} \cite{Hellisto2020}. 

Two of the most common a-priori methods are downsampling and representative periods. Downsampling consists of increasing the coarseness of the time steps used in PSOMs to reduce their size; for example, in SDDP \cite{Pereira1991}, instead of using the available hourly data, weekly averages of demand and inflows are deployed; however, when using coarse time steps researchers lose sight of short-term dynamics; for example, in the case of weekly ones, the daily and hourly patterns of a power system are lost. A case study evaluating this situation is found in \cite{WELSCH2014600}; this means that a model with weekly or daily time steps will not be able to represent the short-term dynamics of the system.

On the other hand, representative periods aim to partition the complete time horizon (e.g. 8760 individual hours of one year) into a weighted average of smaller ones while keeping the coarseness of the time steps (e.g. 7 representative days with 24 individual hours - $7\times24$ different time steps in total). For example, in \cite{Tejada2020}, the authors run 168-hour long models (one week) and use four of these they approximate a complete year, so in this case, we have four representative periods that are used to aggregate the 52 weeks of a year. The main advantage of using representative periods over downsampling is that researchers do not overlook short-term temporal dynamics of the power system; however, the difficulty lies in determining if: first, the chosen periods are representative enough of the original data; and second, when employed in a PSOM, how they translate into accurate results \cite{TEICHGRAEBER2022111984}, i.e., how well they approximate full model results.

However, despite their past usefulness, new trends like the increasing share of variable renewable energy sources (VRES) in power systems, pose a challenge to temporal aggregation techniques as they rely on adding multiple time steps, e.g., daily or weekly averages from hourly data \cite{PFENNINGER20171} \cite{Stenzel2016} \cite{DEANE2014152}, or breaking the inter-temporal linking of time \cite{Battle2013} \cite{GUNER2018310}, but VRES' technical constraints require a highly detailed temporal modeling \cite{COLLINS2017839}. The consequences of these aggregation procedures have been researched \cite{Reichenberg2022} \cite{Pfenninger2014} \cite{Buchholz2020}, and results show that they lead to inaccurate results as they go against these fundamental technical constraints and even those of short-term storage technologies. 

Until now, the inaccuracy of a-priori techniques has not been a matter of concern, as the increasing share of VRES to achieve net-zero power systems has highlighted the limitations of the current temporal aggregation techniques used in PSOMs for planning and policy-making. The limitations of a-priori methods are especially highlighted and quantified in \cite{WOGRIN2023}, where the author illustrates the substantial impact of applying a standard a-priori TSA method (representative periods with k-Means clustering) to approximate a full PSOM. The case study of a simple Economic Dispatch problem shows a 91\% error when the objective function value of the k-Means aggregated model is compared with the full model, while other (a-posteriori) alternatives achieve a 0\% error; this accentuates the need for new and adequate aggregation techniques.

In contrast to a-priori methods, a-posteriori methods do not exclusively rely on PSOM input data but also include information from the PSOM's structure, e.g., the results from partial model runs or relaxations, in the TSA process. A-posteriori methods can be classified as either non-iterative or iterative; non-iterative aim to obtain a fast solution, although with a trade-off between robustness and optimality, they rely on using a relaxed version of the entire model, for example, by reducing or fixing the number of binary variables and then cluster the input data and also identify extreme periods (e.g., those with a high cost in the relaxed model). So, non-iterative methods do not guarantee achieving the optimum but provide \textit{good enough} solutions fast. Iterative methods aim to obtain an aggregated model as close as possible to the complete one with less regard for the computational burden; they work by decomposing the problem, for example, by splitting design and operation or by treating each type of constraint independently and then merging the results. In this way, a-posteriori TSA methods have the potential for significant efficiency gains in terms of solution times while simultaneously maintaining the quality of model results close to the full PSOM (without TSA) \cite{WOGRIN2023},\cite{TEICHGRAEBER2022111984}, \cite{LI2022107697}; this potential has led to an increased research interest about a-posteriori aggregation techniques for PSOM which take into account the structure of the optimization model \cite{Rigo2022} \cite{LI2022107697}\cite{Gonzato2021}.

In advancing a-posteriori methods, \cite{WOGRIN2023} proposes a new approach to temporal aggregation called Basis-Oriented TSA. The Basis-Oriented approach takes a full-hourly solution and splits the input data based on the Simplex basis each time step belongs to; then, for each basis, makes a centroid of the input data and reruns the model using the centroids; this procedure gives an exact approximation to the complete model while retaining computational tractability. Therefore, this approach allows for an aggregated model akin to representative periods, which preserves the objective function value and the solution of the complete model so it has zero error. Another advantage of the Basis-Oriented approach is that, as we demonstrate later in this paper for the case presented in \cite{WOGRIN2023}, it finds an aggregation of minimal cardinality, which means no smaller set of representative periods exists that exactly approximates the model's objective function value. This result is supported by the result in \cite{ZIPKIN1980}, which states that every LP model has a row-aggregation which is equivalent to the complete problem.

In this paper, we extend the research on the Basis-Oriented approach and apply it to significantly more complex instances of PSOMs including additional constraints with both spatial and inter-temporal dependence, thereby bringing the approach closer to real-world applications. In particular, we extend the methodology to network constraints and show that we can directly apply the Basis-Oriented approach to find a zero-error aggregation. Moreover, we extend the approach to constraints with temporal linking, i.e., ramping of thermal generators. However, inter-temporal constraints are more challenging than network constraints, as obtaining a basis is not as straightforward. This complexity arises because the basis duration\footnote{The number of consecutive time periods that cannot be separated without losing model accuracy.} or length (always one hour in \cite{WOGRIN2023}) will vary depending on the system's state (in particular, on inter-temporal dependencies) and it cannot be determined a-priori.

To address this challenge, we develop an algorithm that explores the dual space hour by hour, to identify these bases, which allows us to formulate an aggregated PSOM that is still exact (zero error). In summary, the original contributions of this paper are as follows:
\begin{itemize}
    \item We provide empirical proof that the aggregation found by the Basis-Oriented TSA approach is minimal for the case analyzed in \cite{WOGRIN2023}.
    \item We extend the Basis-Oriented TSA methodology to include network constraints.
    \item We further extend the approach to models with temporal linking (in particular, ramping constraints) and develop an algorithm to find the time-dependent bases to obtain an exact aggregation under such situations.
\end{itemize}

The remainder of this paper is organized as follows: in Section \ref{ssec:agg_repr}, we briefly introduce the original Basis-Oriented TSA approach and show how the aggregation is minimal and unique; then, we extend it to include network and inter-temporal ramping constraints in Section \ref{sec:basisOrientedPSOM}; in Section \ref{sec:results}, we present a case study to validate the TSA methodology and finally, Section \ref{sec:conclusion} presents the conclusions and future research.

\section{Basis-Oriented Time Series Aggregation} 
\label{ssec:agg_repr}
In this section, we first introduce the reader to Basis-Oriented TSA in section \ref{ssec:prevED} by summarizing the main aspects and initial results from \cite{WOGRIN2023}, as it provides the groundwork for our extension to more realistic models. And then, as a novel contribution, in section \ref{ssec:minimBasis} we provide empirical proof that the set of representative periods determined by Basis-Oriented TSA not only obtains an exact aggregation but that it corresponds to the set of minimal cardinality to achieve a zero error, and, that this set is unique. Furthermore, we show that all sets that achieve an exact aggregation with a higher number of representative periods are actually a partition of the original Basis-Oriented result.


\subsection{Previous Results: Economic Dispatch}
\label{ssec:prevED}

Consider a stylized PSOM for the Economic Dispatch (ED) problem where the objective is to minimize the total system operation cost while satisfying demand and scheduling the generation units within their operational bounds. The set of constraints in \eqref{eqn:ED_Full} corresponds to a full model run of the ED problem for every hour $k$ ($k=1,2, \dots, 8760$). Full model results are optimal production decisions $p_{g,k}$ for each hour $k$ and each generator $g$.
In contrast, the set of constraints in \eqref{eqn:ED_Agg} represents an aggregated PSOM. Instead of $k$, this formulation uses representative periods $r$, a subset of $k$. 
To avoid the computational complexity or potential intractability of running the full model, the modeler must choose $r$ such that $|r|\ll |k|$ while obtaining results similar to those of the full model. In order to determine $r$, TSA methods are applied. Common results of a TSA method for a specific number of $r$ are: representative values for the demand $D_r$; a factor $W_r$, which corresponds to the weight (or a number of occurrences) of each $r$ so that the aggregated model represents the ED problem for a full year; and other data $\overline{P}_{g,r}$.


  \begin{subequations}\label{eqn:ED_Full}
  \begin{align}
      min    \sum_{g,k} C_g p_{g,k} + \sum_k C^{nsp} nsp_k \label{eqn:ED_FullObj}\\
  \text{s.t. } \sum_{g} p_{g,k} + nsp_k = D_k \quad \forall k \label{eqn:ED_FullBalance} \\
  \underline{P}_{g} \leq p_{g,k} \leq \overline{P}_{g,k}  \quad \forall g,k \label{eqn:ED_FullBounds}
    \end{align}
    \end{subequations}

    \begin{subequations}\label{eqn:ED_Agg}
    \begin{align}
      min    \sum_{g,r} C_g p_{g,r} + \sum_r C^{nsp} nsp_r  W_r \label{eqn:ED_AggrObj}\\
  \text{s.t. } \sum_{g} p_{g,r}+nsp_r = D_r \quad \forall r \label{eqn:ED_AggrBalance} \\
  \underline{P}_{g} \leq p_{g,r} \leq \overline{P}_{g,r} \quad \forall g,r \label{eqn:ED_AggrBounds}
    \end{align}
    \end{subequations}


In \cite{WOGRIN2023}, the author shows that the full hourly ED model, which uses the whole time series of 8760 individual hours and a simple generation mix (one thermal, one wind plant, and the option of non-supplied energy), can be approximated exactly using only three representative hours. Instead of finding these three representative hours by the proximity of the data points (such as by k-means), they are determined by the Basis-Oriented approach - identifying the Simplex basis to which each hour belongs to\footnote{Since there are no inter-temporal constraints in this model, every hour is assessed separately. Hence, we can speak of an hourly basis. In the case study presented, there were only 3 out of 8760 possible different hourly bases. The 3 representative hours were then obtained as the averages of all hours that belonged to the same Simplex basis.}. However, \cite{WOGRIN2023} does not address the following questions: Are there other aggregations with three clusters that also yield zero error, and if so, what characterizes them? In the following section, we attempt to answer this question and further characterize the results of Basis-Oriented TSA.

\subsection{Minimal Number of Bases and Uniqueness}
\label{ssec:minimBasis}

As previously mentioned, Basis-Oriented TSA achieves an exact approximation of a full hourly ED, using three representative hours only; despite their work ensuring zero error, it gives no insight concerning the existence and number of aggregations when using a different number of clusters. So, this section aims to assess whether this is the only exact aggregation with zero error using three clusters and how other exact aggregations, with a different number of clusters\footnote{It is important to note that we use the word clusters, particularly in this section, in the word's general meaning and not the one specific to the machine learning community.}, relate to the one determined under Basis-Oriented TSA. We do so employing exhaustive numerical enumeration. 

For this experiment and for the sake of simplicity, we limit the total time horizon of the ED problem to 12 hours. Hence, the hourly ED problem in \eqref{eqn:ED_Full} has $|k|$ = 12 time steps. This yields the exact solution to the ED problem and a reference to compare the quality of the approximation under the aggregated model in \eqref{eqn:ED_Agg}. 
For the aggregated model, there exists a plethora of combinations to group the original data into 1$\leq|r|\leq$12 clusters. For example, there is only one possible combination to group the data into one cluster. However, according to the Stirling Numbers of the Second Kind \cite{Graham1994}, there are 2047 combinations of organizing the 12 data points into two clusters, and so on until 12 clusters. Table \ref{tab_minUnique} presents the relation between the number of clusters and the number of possible combinations; in total there are \mbox{4 213 597} different ways of clustering the 12 hourly data points. We run the aggregated model for every possible combination and study the quality of the solution in terms of the error in the objective function. Table \ref{tab_minUnique} shows the results of those runs.
Therein, column \textit{Clusterings with No Error} corresponds to the number of combinations that exactly approximated the full hourly model. For example, from the 2047 possible combinations to separate 12 hours into two clusters, none yielded a 0\% error when used in an aggregated model. From the 8526 ways to organize the original data into 3 clusters, only one yields an exact approximation of the original hourly model, i.e., the one obtained through basis-oriented TSA. This proves that basis-oriented TSA yielded the best possible aggregation for this data and the economic dispatch problem (as there is no exact aggregation with lower cardinality), and it is unique for this number of clusters.  

\begin{figure}[htbp]
     \centering
     \begin{subfigure}[b]{0.23\textwidth}
         \centering
         \includegraphics[width=\textwidth]{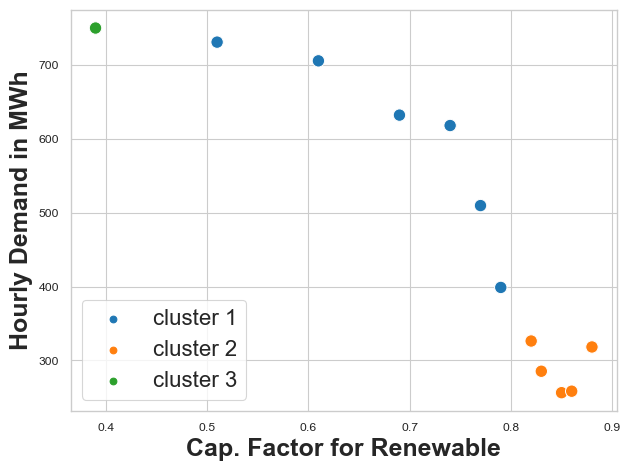}
         \caption{Three clusters with no error}
         \label{fig_three_clusters}
     \end{subfigure}
     \hfill
     \begin{subfigure}[b]{0.23\textwidth}
         \centering
         \includegraphics[width=\textwidth]{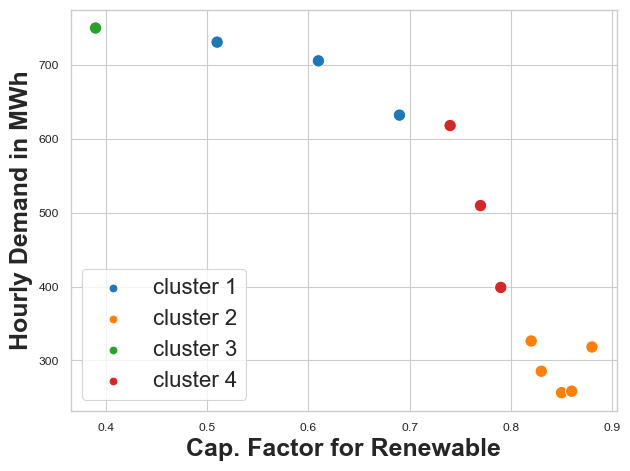}
         \caption{Four clusters with no error}
         \label{fig_four_clusters}
     \end{subfigure}
    \label{fig_basis}
    \caption{Two examples of clusterings with no error: (a) with 3 clusters and (b) with 4 clusters}
\end{figure}

Now, let us analyze all the other combinations of aggregations with a cardinality higher than three that have also led to an exact aggregation. In Fig. \ref{fig_three_clusters} the Basis-Oriented aggregation is presented and in Fig. \ref{fig_four_clusters} we see a four clusters aggregation which also yields zero error, but clusters one and four in Fig. \ref{fig_four_clusters} are just a division of cluster two from Fig. \ref{fig_three_clusters}. This situation occurs for all other aggregation with zero error and cardinality higher than three: all of these aggregations correspond to sub-divisions of the basis-oriented 3-cluster separation. This result is quite remarkable because it shows that not only does the Basis-Oriented approach yield the best possible approximation but that other exact aggregations come from the Basis-Oriented one.

\begin{table}[!ht]
    \caption{Combinatorial Aggregation Results}
    \label{tab_minUnique}
    \centering
    \begin{tabular}{|c|c|c|}
    \hline
        \textbf{Clusters $|r|$} & \textbf{Possible Clusterings} & \textbf{Clusterings with No Error} \\ \hline
        1 & 1 & 0 \\ \hline
        2 & 2 047 & 0 \\ \hline
        3 & 8 526 & 1 \\ \hline
        4 & 611 501 & 256 \\ \hline
        5 & 1 379 400 & 3280 \\ \hline
        6 & 1 323 652 & 10795 \\ \hline
        7 & 627 396 & 14721 \\ \hline
        8 & 159 027 & 9597 \\ \hline
        9 & 22 275 & 3108 \\ \hline
        10 & 1 705 & 498 \\ \hline
        11 & 66 & 37 \\ \hline
        12 & 1 & 1 \\ \hline
    \end{tabular}
\end{table}

Based on these results, which underline the potential of Basis-Oriented TSA, we set out to extend the methodology to more realistic models and expand this work by adding time-linking ramping constraints and network constraints, which are commonly found in energy system models and discuss the arising challenges.

\section{Extension of Basis-Oriented TSA framework}
\label{sec:basisOrientedPSOM}

Considering that PSOM vary significantly in their formulation depending on the aspects being analyzed, our goal in this section is to frame the varieties of constraints and technical aspects we consider in our extension to the Basis-Oriented approach. 
Commonplace formulations of PSOMs include network aspects, e.g., the DC Optimal Power Flow (OPF) \cite{Carpentier1979} or the AC OPF or its convexification \cite{Low2014}, which are widely used by operators for day-ahead scheduling and market clearing; or security-constrained unit commitment (SCUC) \cite{Fu2005} including inter-temporal constraints. Therefore, this paper focuses on extending the ED formulation and the Basis-Oriented TSA to: first, consider network flow constraints in section \ref{ssec:basisNetwork}, and then ramping constraints in section \ref{ssec:ramping}.

\subsection{Basis-Oriented TSA: Network Constraints}
\label{ssec:basisNetwork}

In this section, we extend the Basis-Oriented aggregation to include network flow constraints and show that it still finds an exact aggregation of the input-data. The full hourly model for the ED problem with network flows is presented in \eqref{eqn:nwFull}. The additional term in the objective function compared to \eqref{eqn:ED_Full} corresponds to transmission costs while constraints \eqref{eqn:nwFull_Balance}-\eqref{eqn:nwFull_BoundsOut} represent the nodal equations, the import limits and the export limits respectively. To find a Basis-Oriented aggregated model, we group the hours that have the same hourly basis and then form the centroid from each of these groups; the only difference this time is the total number of bases under the additional constraints. Because of the proof in \cite{WOGRIN2023}, and in the absence of temporal linking constraints, this aggregation has zero error if compared to the full model.
  \begin{subequations}
    \label{eqn:nwFull}
    \begin{align}
        \min \sum_{k,g} C^{g} p_{g,k}  + \sum_{k,i} C^{nsp} nsp_{i,k} + \sum_{k,i,j}C^Nf_{k, j, i}  \label{eqn:nwFull_Obj}\\
        \text{s.t. \;} \underline{P_{w}} \leq p_{w,k} \leq CF_{k} \overline{P_{w}} \quad \forall{k,w} \label{eqn:nwFull_BoundsW}\\
        \underline{P_{t}} \leq p_{t,k} \leq \overline{P_{t}} \qquad \quad \forall{k,t}  \label{eqn:nwFull_BoundsT}\\
        \sum_{j} f_{k, j, i} - \sum_{j} f_{k, i, j} + nsp_{i,k} + \sum_{g \in i} p_{g,k} = D_{k,i} \forall{k,i} \label{eqn:nwFull_Balance}\\
        \sum_{j} f_{k, i, j} \leq \overline{F_{l}} \quad \forall{k,i} \label{eqn:nwFull_BoundsOut}\\
        \sum_{j} f_{k, j, i} \leq \overline{F_{l}} \quad \forall{k,i} \label{eqn:nwFull_BoundsIn}
    \end{align}
   \end{subequations}

\subsection{Basis-Oriented TSA: Ramping Constraints}
\label{ssec:ramping}

Especially in PSOM, chronology matters, and this is why it is so important to account for it in TSA methods. A common mathematical constraint found in PSOMs and that links together multiple temporal periods is a ramping constraint. This is used to represent the technical characteristic of thermal generators. In the case of ramping constraints, the temporal linking arises because of the maximum changes in production a thermal generator $t$ can tolerate between two consecutive periods; for example, if $p_{t,k}$ corresponds to the power produced by generator $t$ during hour $k$, a ramping up constraint would look like $p_{t,k} - p_{t,k-1} \leq RU$  while a ramping down constraint would be $p_{t,k-1} - p_{t,k} \leq RD$ where $RU$ and $RD$ are parameters of the generator.

To analyze the impact temporal linking has on the Basis-Oriented approach, we use the model in \eqref{eqn:rmpFull}, which corresponds to the model used in \cite{WOGRIN2023} with added ramping constraints for the thermal generator $t$.
  \begin{subequations}\label{eqn:rmpFull}
  \begin{align}
      min    \sum_{g,k} &C_g p_{g,k}  + \sum_k C^{nsp} nsp_k\label{eqn:rmpFull_Obj}\\
  \text{s.t. } \sum_{g} p_{g,k} + nsp_k &= D_k & \forall k \label{eqn:rmpFull_Balance} \\
    p_{t,k} - p_{t,k-1} &\leq RU  & \forall{k,t} \label{eqn:rmpFull_RU}\\
    p_{t,k-1} - p_{t,k} &\leq RD  & \forall{k,t} \label{eqn:rmpFull_RD}\\
  \underline{P}_{g} \leq p_{g,k} & \leq \overline{P}_{g,k}  & \forall g,k \label{eqn:rmpFull_Bounds}
    \end{align}
    \end{subequations}

Temporal linking constraints pose a challenge to the application of the Basis-Oriented TSA as it is no longer possible to assess the data of different hours separately - because due to active ramping constraints in the optimization model, the hourly periods are no longer independent. In Basis-Oriented TSA we therefore need to consider that a basis $b$ might consist of multiple consecutive periods; in other words, some representative periods might be longer than one hour (because of active ramping constraints). Identifying those periods and their corresponding duration is not trivial a-priori.

To give an example, if a ramping constraint is active in a given time step, it will bind variables from two different time steps. Another way of thinking about this is considering the structure of the ramping constraint. Without it, the thermal generation in one time step is independent of the others, that is, a line with slope zero in the $p_{t,k}, p_{t,k-1}$ plane; but with ramping, they now represent a line with a slope of 1. This extends to the number of \textit{joined} periods the ramping constraint is active. From this, we obtain a \textit{basis duration}, which corresponds to the set of consecutive periods that depend on each other at the optimum and therefore cannot be split. 

This leads to another contribution of this paper - how to identify these joined bases and their corresponding basis durations, which we describe in the following section.

\subsection{Basis-Oriented TSA: Basis Identification and Heuristic Algorithm}
\label{ssec:rampingalg}

In this section, we discuss how to identify bases for PSOMs with temporal linking constraints, such as ramping, and propose a heuristic algorithm to that purpose. This heuristic algorithm uses the dual solution of the full model to identify bases. We are aware that in real-life examples, this dual solution is not available. However, our goal here is to demonstrate that even for temporal linking and flow limit constraints the Basis-Oriented approach still finds an aggregation with zero error with respect to the objective function value. Moreover, using dual information, such as marginal costs\footnote{The marginal cost is the dual variable of the balance equation and corresponds to the cost of serving one additional unit of demand in that hour.} is very easy to understand and provides great intuition. In future research, we plan to develop a methodology that allows us to identify these bases using input data only.

To account for temporal linking, we use the marginal cost as a proxy to identify bases. In the simple case of \cite{WOGRIN2023}, where no temporal interlinking constraints were present, only three different marginal costs could be observed, i.e., the variable cost of the wind generator, the thermal generator, and the cost of non-supplied energy. However, introducing ramping constraints complicates this and causes a plethora of additional marginal costs to appear, which do not correspond to any of the variable costs of the system's generation units. This is because, with active ramping constraints, the marginal cost of one hour is affected by the system operation of preceding or posterior hours.

For example, consider a four-hour dispatch of a single-node system with a wind turbine and a thermal generator with variable costs 3~\texteuro/MWh and 24~\texteuro/MWh respectively. Table \ref{tab_mc_example:no_ramp} shows system demand, the cost-optimal productions of the generators, and the arising marginal costs. Note that between hours 3 and 4, the thermal plant has a production increase that exceeds 100 MW. Let us repeat this example but with an additional ramp-up constraint with a limit of 100 MW for the thermal generator. In Table \ref{tab_mc_example}, we now observe that production and marginal costs have changed. Since the ramping constraint prohibited that large jump from hour 3 to 4, the thermal plant had to increase its production earlier on - replacing the cheaper wind production in hours 2 and 3. The marginal cost of 66~\texteuro/MWh occurs because the thermal generator has to increase its production two hours before (indicated in green in Table \ref{tab_mc_example}), plus the cost of serving the additional demand (in blue), and minus the cost of the displaced production from wind (in red), so $\textcolor{green}{2\cdot24}+\textcolor{blue}{1\cdot24}-\textcolor{red}{2\cdot3} = 66$.

This example illustrates how temporal interlinking directly affects the marginal cost of the system not only in the period where it occurs but also in the production in periods linked to it. Generalizing the previous case, we find that for our model, the marginal cost of any given hour is a linear combination of the variable costs of the generating units with integer coefficients $a,b,c$ because of the absence of losses.

\begin{table}[htbp]
    \caption{Example: No Ramping Constraint}
    \label{tab_mc_example:no_ramp}
    \centering
    \begin{tabular}{|c||c|c|c|c|c|}\hline
    \textbf{Period}  & \textbf{Demand} & \textbf{Thermal}   & \textbf{Wind} &  \textbf{Ramp-Up}   & \textbf{MC} \\\hline
    $1$	 & 170.2	 & 130.2	 & 40.0    &  -	     & 24	 \\\hline
    $2$	 & 176.0  & 171.0 & 5.0    & 40.8	 & 24	 \\\hline
    $3$	 & 281.7 & 266.0  & 15.7 & 100.0	 & 24 \\\hline
    $4$	 & 391.0 & 371.0	 & 20.0    & 105.0	 & 24	 \\\hline
    \end{tabular}
\end{table}

\begin{table}[htbp]
    \caption{Example: With Ramping (blue/green/red indicate production change and relation to marginal cost)}
    \label{tab_mc_example}
    \centering
    \begin{tabular}{|c||c|c|c|c|c|}\hline
    \textbf{Period}  & \textbf{Demand} & \textbf{Thermal}   & \textbf{Wind} &  \textbf{Ramp-Up}   & \textbf{MC} \\\hline
    $1$	 & 170.2	 & 130.2	 & 40.0    &  -	     & 24	 \\\hline
    $2$	 & 176.0  & 171.0 \textcolor{green}{(+1)} & 5.0 \textcolor{red}{(-1)}    & 40.8	 & 24	 \\\hline
    $3$	 & 281.7 & 271.0 \textcolor{green}{(+1)} & 10.7 \textcolor{red}{(-1)}    & 100.0	 & 45	 \\\hline
    $4$	 & 391.0 \textcolor{blue}{(+1)} & 371.0	 & 20.0    & 100.0	 & 66	 \\\hline
    \end{tabular}
\end{table}

\begin{equation}
    MC_{k} = a VC_{nsp} + b VC_{t} + c VC_{w} \quad : \quad a,b,c \in \mathbf{Z}
    \label{eq_marg_cost}
\end{equation}

\begin{figure}[h]
\centering
\includegraphics[width=0.5\textwidth]{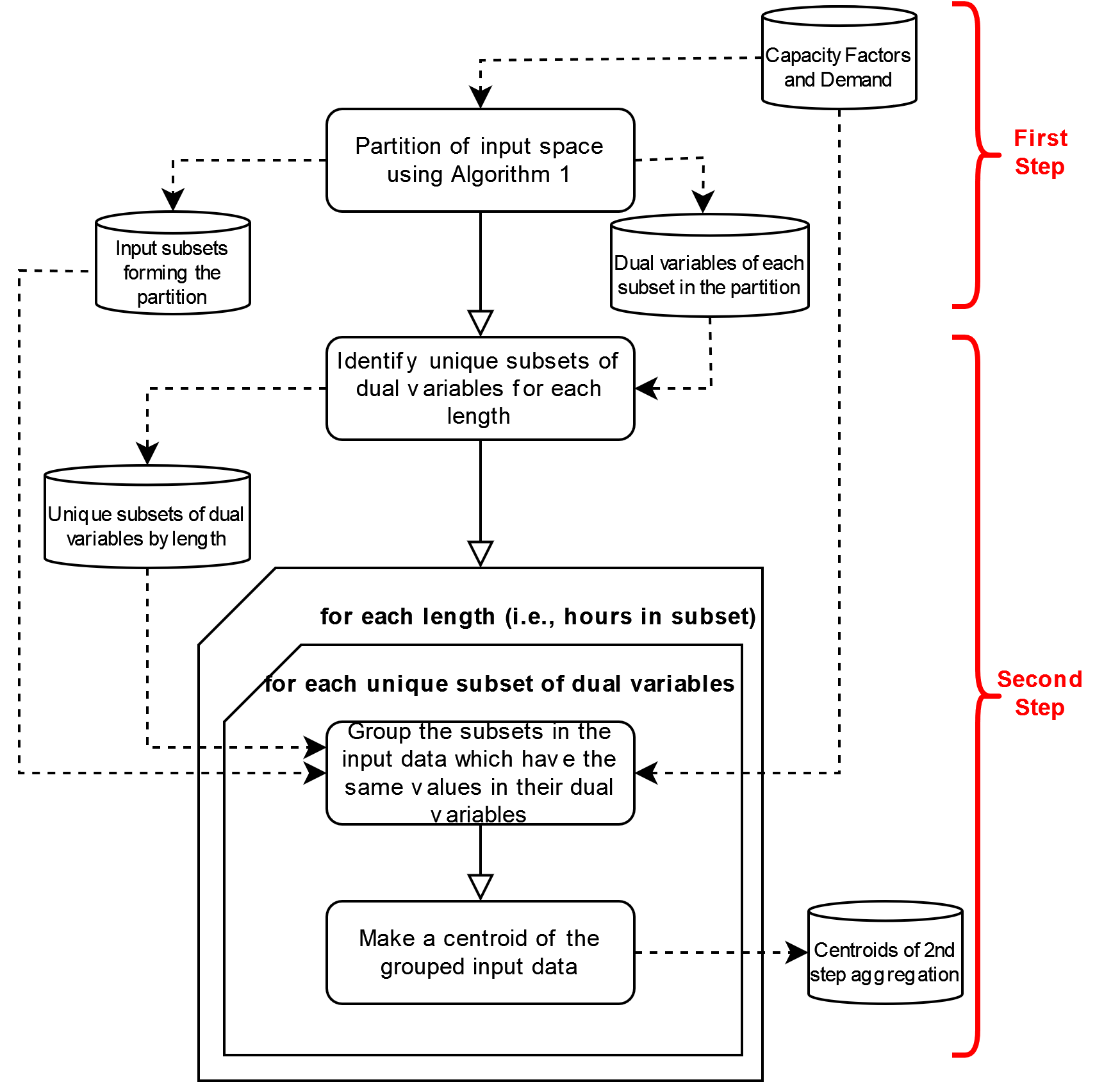}
\caption{Data-aggregation process with temporal linking constraints}
\label{fig_agg_steps}
\end{figure}

With this observation about the marginal cost in mind, we develop a heuristic algorithm, i.e., Algorithm~\ref{alg:cap}, detailed in the Appendix. Algorithm~\ref{alg:cap} partitions the time horizon into groups of hours that cannot be separated (because of active ramping constraints) using information from the duals (as demonstrated in the example). The arising partition of the input data space - based on these groups of hours - will then be used to build an aggregated model using dual variables. Fig. \ref{fig_agg_steps} describes the complete data-aggregation process that we employ with temporal linking constraints.

The \textbf{First Step} of the data-aggregation process consists of loading the input data, capacity factors, and demand and running Algorithm \ref{alg:cap}\footnote{Algorithm \ref{alg:cap} looks for a time period where the system's marginal cost does not correspond to the variable cost of any of the generation resources. Looking at the dual variables of the ramping constraints in the periods after and before, we identify if the marginal cost comes from an up or downward ramping constraint. The length of the ramp is determined by the closest integer multiple to the variable cost of the thermal generator; this is based on the empirical work in \cite{DiTondo23} where the authors found that, in the case of marginal costs higher than $VC_t$ but lower than $VC_{nsp}$, the integer multiple of $VC_t$ corresponded to the number of periods in which one of the ramping constraints was binding.}; from this, we obtain a temporal partition of the inputs and the dual variables with subsets of different lengths\footnote{During our experiments, we found that the more constrained a model is, for example, a low ramping capacity, the higher the number of hours in each element.} (measured in numbers of consecutive time periods i.e. hours); then, we want to identify the bases; hence, for each length (e.g. for all 3-hour partitions) we check which of those partitions have the exact same dual variables, and how many different \textit{basis} (i.e., subsets that belong to different partitions according to Algorithm~\ref{alg:cap} but have the same values in the dual space) there are. For example, in the 8760 hour time horizon, there might be 189 total subsets with a 3-hour duration. However, when comparing their dual variables, there are only 7 different subsets of dual variables. Hence, there are only 7 bases of length 3 hours.

After that, in the \textbf{Second Step}, we want to aggregate the data belonging to the same basis. Hence, we group the input data using the unique subsets of dual variables found in the first step and make a centroid for each of them with a weight corresponding to the number of elements in the subset, each of these centroids is a representative period in our aggregated model.

In this subsection, we have presented a data aggregation process that decreases the computational complexity of the model with time-linking constraints and performs a temporal aggregation while maintaining a zero error when compared to the complete model. 


\subsection{Basis-Oriented TSA: Network and Ramping Constraints}
\label{ssec:flowLimRamping}

In this section, we consider the model in \eqref{eqn:rmpNwFull}, which includes flow-limit and ramping constraints. As previously illustrated, flow limits are easier to handle than ramping because there is no temporal interlinking, and they operate in a way analogous to hourly available capacity constraints. So we can apply Algorithm \ref{alg:cap} making a small tweak: we now have to look at the marginal costs of every bus when looking for the periods that must go together, however, the aggregation procedure is the same as in Fig. \ref{fig_agg_steps}.

\allowdisplaybreaks

  \begin{subequations}
    \label{eqn:rmpNwFull}
    \begin{align}
        \min \sum_{k,g} C^{g} p_{g,k}  + \sum_{k,i} C^{nsp} nsp_{i,k} + \sum_{k,i,j}C^Nf_{k, j, i}   \label{eqn:rmpNwFull_Obj}\\
        \text{s.t. \;} \underline{P_{w}} \leq p_{w,k,i} \leq CF_{k} \overline{P_{w}} \quad \forall{k,w,i} \label{eqn:rmpNwFull_BoundsW}\\
        \underline{P_{t}} \leq p_{t,k,i} \leq \overline{P_{t}} \qquad \quad \forall{k,t,i}  \label{eqn:rmpNwFull_BoundsT}\\
        \sum_{j} f_{k, j, i} - \sum_{j} f_{k, i, j} + nsp_{i,k} +\sum_{g \in i} p_{g,k} = D_{k,i} \forall{k,i} \label{eqn:rmpNwFull_Balance}\\
        p_{t,k} - p_{t,k-1} \leq RU  \quad \forall{k,t} \label{eqn:rmpNwFull_RU}\\
        p_{t,k-1} - p_{t,k} \leq RD  \quad \forall{k,t} \label{eqn:rmpNwFull_RD}\\
        \sum_{j} f_{k, i, j} \leq \overline{F_{l}} \quad \forall{k,i} \label{eqn:rmpNwFull_BoundsOut}\\
        \sum_{j} f_{k, j, i} \leq \overline{F_{l}} \quad \forall{k,i} \label{eqn:rmpNwFull_BoundsIn}
    \end{align}
   \end{subequations}

\section{Numerical Results}
\label{sec:results}

This section contains experimental results obtained from applying the extended Basis-Oriented approach to the 3-bus system described in Fig. \ref{fig_three_buses_diagram}. The input data and parametrization can be found in \cite{BasisGithub}; the system has a maximum demand of 1000 MW while the generation units have the parametrization presented in Table \ref{tab_gen_units} and non-supplied power has a cost of \texteuro/MWh 5000. The data are such that the demand could be supplied entirely with the thermal generator; also, the limits in the lines represent the difficulty of getting renewable energy from the wind turbine into the grid as is usually the case in real-world systems.

\begin{table}[htbp]
    \caption{Parameters of generating units}
    \label{tab_gen_units}
    \centering
    \begin{tabular}{|c|c|c|}\hline
    \textbf{\quad}      & \textbf{Wind Turbine}  & \textbf{ Thermal Generator}\\\hline
    \textbf{Installed Capacity (MW)} & 500 & 1000 \\ \hline
    \textbf{Variable Cost (\texteuro/MWh)} & 3 & 24 \\ \hline
    \textbf{Ramp limit (MW/h)} & - & 100 \\ \hline
    \end{tabular}
\end{table}

\subsection{Results with Network Constraints}
\label{ssec:resBasisNetwork}

\begin{figure}[h]
\centering
\includegraphics[width=0.4\textwidth]{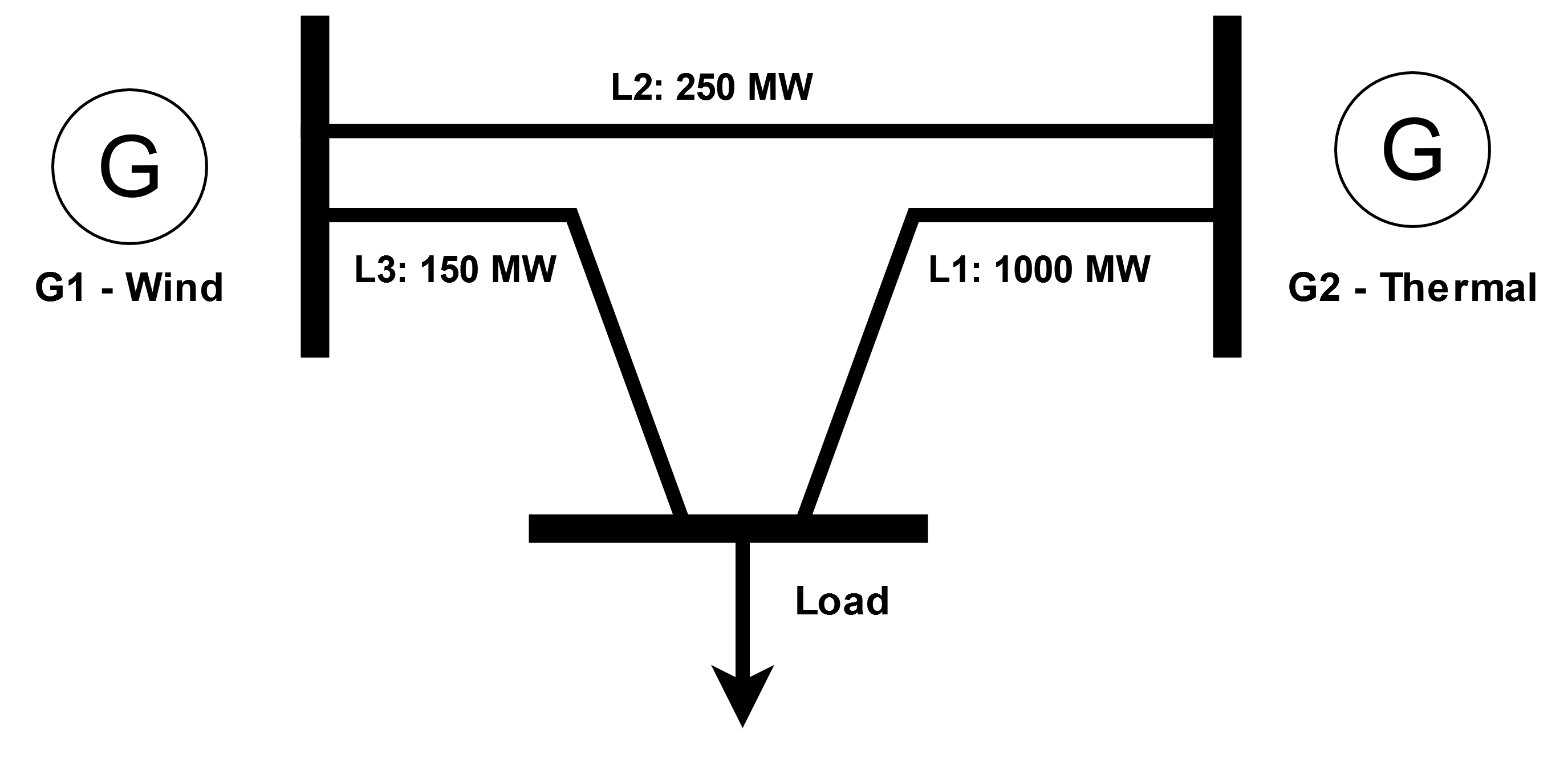}
\caption{Three-Bus system diagram with line limits (MW) and location of generators}
\label{fig_three_buses_diagram}
\end{figure}

We applied the Basis-Oriented TSA to the system presented in Fig. \ref{fig_three_buses_diagram} whose data and parameter values are in \cite{BasisGithub}. For simplicity, we consider that load is present only in one node of the system while the other two are net exporters. The model is solved for one year with hourly time steps to obtain a benchmark against which the temporally aggregated model is compared to\footnote{The full hourly model poses no tractability challenges for small test systems.}. After applying the Basis-Oriented procedure to the network constraints case, we obtained the same results from the temporally aggregated model and the full-hourly one; this means that Basis-Oriented TSA yields an exact aggregation and achieved a reduction in the number of required hours of 1747 times (8736/5) because we only required five representative hours to exactly represent the whole year, as shown in Table \ref{tab_3BNR}. This corresponds to a reduction of the number of variables of the aggregated model of 3 orders of magnitude. The relationship between the computation time and the number of hours is not linear, however, it is one of the determining factors in solving such models in a feasible amount of time.

\begin{table}[htbp]
\caption{Three-Bus System}
\begin{center}
\begin{tabular}{|c|c|c|c|}
\cline{1-4}
\textbf{}& \makecell{\textbf{Aggregated}\\ \textbf{Model}} & \makecell{\textbf{Complete} \\ \textbf{Model}}  & \makecell{\textbf{Size} \\ \textbf{Reduction (\%)}}\\ 
\cline{1-4}
\makecell{Number of  Variables}& 45 & 78 624 & \multirow{2}{*}{99,94\%}\\
\cline{1-3}
\makecell{Number of  Constraints}& 55 & 96 096 & \\
\cline{1-4}
\end{tabular}
\label{tab_3BNR}
\end{center}
\end{table}

\begin{figure}[h]
\centering
\includegraphics[width=0.4\textwidth]{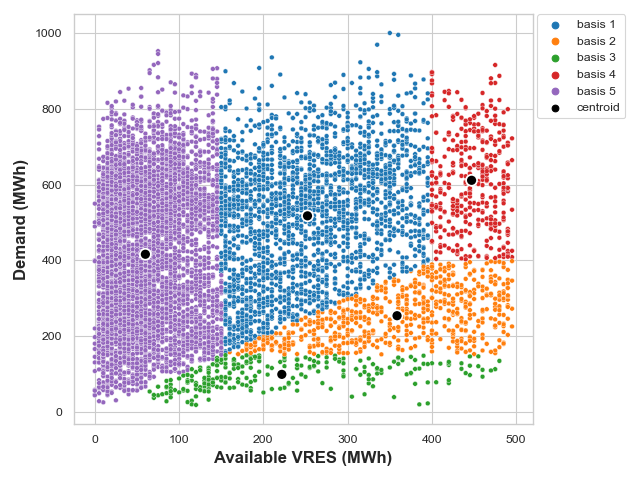}
\caption{Bases (basis 1-5) in the input-data space - Three bus system}
\label{fig_three_buses_basis}
\end{figure}

To analyze the results obtained from the Basis-Oriented approach, we plot the resulting bases in the input-data space in Fig. \ref{fig_three_buses_basis} (renewable energy availability on the x-axis and demand on the y-axis). Comparing the results with those in \cite{WOGRIN2023} (where we only had 3 linearly separable bases), we see that in the three-bus network case, the bases are not linearly separable; this stresses the point that \textit{a-priori} approaches for representative periods selection are insufficient in PSOM. These results also strengthen the argument that a partitioning of the inputs from a minimum distance perspective (i.e.: k-Means or k-Medoids clustering) is insufficient to approximate the complete model effectively. 

In the network case, we obtain five bases because of the additional\footnote{Compared to the case in \cite{WOGRIN2023}.} constraints in the model; in this case, the additional bases come from the hours when lines one or two are congested, which means that one or two of the constraints in \eqref{eqn:rmpNwFull_BoundsOut} \eqref{eqn:rmpNwFull_BoundsIn} are binding. For this case, it is worth noting that each of the bases corresponds to an operational situation in the power system: \textit{basis 1} and \textit{basis 5} are hours where the demand is higher than the available renewable energy; for \textit{basis 1} those hours are the ones where no line is congested while \textit{basis 5} are hours where L1 is at its limit, so a part of the renewable production is routed through L2 and L3; an analogous situation happens with \textit{basis 2} and \textit{basis 3}, but in this case, the demand is lower than the available renewable energy so \eqref{eqn:rmpNwFull_BoundsW} is not binding; finally, \textit{basis 5} correspond to hours where both L1 and L2 are at their limit. In this situation, it makes no difference whether or not \eqref{eqn:rmpNwFull_BoundsW} is binding for the wind turbine and that's why \textit{basis 5} includes hours with demand higher and lower than the available renewable energy. Of course, depending on the input data, some constraints might never be binding at all (e.g., lines have much capacity when compared to the generators they serve, like L3 and the thermal generator in our test case); that is why to precisely approximate these models, both the input data and the model's structure must be considered.

\subsection{Results with Ramping Constraints}
\label{ssec:resRamp}

Running the full hourly model \eqref{eqn:rmpFull} (single-node problem with ramping constraints) with the parametrization in \cite{BasisGithub}, results in the marginal costs shown in Table \ref{tab_decomp}, which also includes their integer coefficients ($a,b,c$) corresponding to each of the generation units' variable cost and the non-supplied power penalty. Now we want to generate an aggregated model and present two different cases: Case A, where we ignore the fact that there is temporal linking in the TSA process; and, Case B, where we account for temporal linking constraints employing the heuristic data-aggregation process presented in Section \ref{ssec:rampingalg}.

\begin{table}[htbp]
    \caption{Possible Marginal Costs, frequency of occurrence and integer coefficients ($a,b,c$) of variable costs}
    \centering
    \begin{tabular}{|c|c|c|c|c|}
    \hline
        \textbf{MC Value} & \textbf{Frequency} & $a$ & $b$ & $c$ \\ \hline
        3 & 1868 & 0 & 0 & 1 \\ \hline
        24 & 6203 & 0 & 1 & 0 \\ \hline
        45 & 247 & 0 & 2 & -1 \\ \hline
        66 & 137 & 0 & 3 & -2 \\ \hline
        87 & 81 & 0 & 4 & -3 \\ \hline
        108 & 24 & 0 & 5 & -4 \\ \hline
        129 & 7 & 0 & 6 & -5 \\ \hline
        5000 & 99 & 1 & 0 & 0 \\ \hline
    \end{tabular}
    \label{tab_decomp}
\end{table}


\subsubsection{Case A (ignoring temporal linking in TSA)}
\label{sssec:CaseA}

We apply the original Basis-Oriented TSA approach (which implicitly assumes hourly bases) from \cite{WOGRIN2023} to the single-node system in \eqref{eqn:rmpFull} and obtain 53 unique one-hour-long combinations in the dual space. So, 53 different hourly bases. Then we partition and aggregate the input data within these 53 hourly bases, which corresponds to an aggregation of 99,39\% with respect to the number of variables; and run the aggregated models\footnote{As opposed to the TSA procedure, the actual models include ramping constraints.} (using cluster centroids and weights) obtaining an objective function value that is 22\% lower than the full model; a summary of the results is presented in Table \ref{tab_SNR}; this shows that assuming a fixed one-hour duration/length of the basis (as done in \cite{WOGRIN2023}) is insufficient to obtain an exact solution under temporal linking constraints.

\begin{table}[htbp]
\caption{Case A: Comparison full versus aggregated model}
\begin{center}
\begin{tabular}{|c|c|c|}
\hline
\textbf{Variable}& \textbf{Aggregated Model} & \textbf{Complete Model}\\ \hline
        Objective Function (\texteuro) &  61 722 846  &  75 269 408  \\ \hline
        Thermal Generation (MWh) &  2 403 380  &  2 537 604  \\ \hline
        Wind Generation (MWh) &  1 347 240  &  1 210 869  \\ \hline
\end{tabular}
\label{tab_SNR}
\end{center}
\end{table}

\subsubsection{Case B (accounting for temporal linking in TSA)}
\label{sssec:CaseB}

Given Case~A, in the following we present an extension of the original Basis-Oriented TSA, which is an original contribution of this paper. The general idea is to drop the assumption of a fixed one-hour basis length and allow for a variable number of hours in each basis. 
The number of hours that we must group together comes from the heuristic process presented in Section \ref{ssec:rampingalg}, where we showed that the system's marginal cost has information concerning how many, and which hours should be kept together in an aggregated model. Column $b$\footnote{The integer coefficient of the thermal generator's variable cost in the decomposition of the marginal cost.} in Table \ref{tab_decomp} indicates how many hours a ramp-up (ramp-down) constraint was active before (after) a given hour. This means that the bases will not be of the same length since it depends on the heuristic aggregation carried out by Algorithm~\ref{alg:cap}. After applying Algorithm \ref{alg:cap} (1st step of aggregation process presented in Fig.~\ref{fig_agg_steps}) we obtain a partition of the input data into 489 subsets with lengths ranging from 1 to 53 hours. The results are summarized in Table \ref{tab_agg1}. For example, in the 8736-hour time horizon, there are 189 subsets with a 3-hour length, 6390 subsets with a 1-hour length etc. The 1-hour long subsets represent situations where only the balance constraint is binding (and ramping constraints are inactive). The column \textit{Obj. Func. Avg.} refers to the average objective function value of all of the subsets of a given length. For example, on average a 3-hour subset yields an objective function value of 22 624\texteuro.  

We include the results in Table~\ref{tab_agg1} because they allow for interesting observations. For example, they indicate the impact of extreme periods on the objective function value as some of the subsets with the highest objective function values have the lowest frequency (e.g., 43 hours which only appears once). Also, note that the lengths do not correspond to the typically used aggregations of days or weeks (e.g., 24 and 168 hours, respectively), which might indicate that using typical representative days in aggregated models is not the most efficient way of aggregation. Moreover, a higher length (number of hours within a subset) does not necessarily imply a higher objective function cost (e.g., the 43-hour subsets have a higher average cost than the 53-hour subsets).

\begin{table}[htbp]
    \caption{Results data-aggregation process: partition subsets (1st Step) and bases (2nd step)}
    \begin{center}
    \begin{tabular}{|c|c|c|c|}
    \hline
        \textbf{Length}         & \textbf{\# of Subsets} & \textbf{\# of Bases} &  \textbf{Obj. Fun. Avg.} \\ 
        \textbf{(\# of hours)}  & \textbf{1st Step}      & \textbf{2nd Step}    &  \textbf{1st Step}  \\        
        \hline
        1 & 6390 &  2   &4 526 \\ \hline
        2 & 1 &   1  &34 632 \\ \hline
        3 & 189 & 7  &22 624 \\ \hline
        4 & 133 & 9    &42 147 \\ \hline
        5 & 77 &  10   &44 791 \\ \hline
        6 & 41 &  8   &42 718 \\ \hline
        7 & 12 &  4   &59 825 \\ \hline
        8 & 3 &   3  &196 191 \\ \hline
        10 & 6 &  6   &173 078 \\ \hline
        11 & 2 &  2   &215 971 \\ \hline
        12 & 10 &  7   &312 782 \\ \hline
        13 & 2 &  2   &496 481 \\ \hline
        14 & 4 &  4   &349 656 \\ \hline
        16 & 1 &  1   &363 219 \\ \hline
        21 & 1 &  1   &186 974 \\ \hline
        27 & 2 &  2   &363 894 \\ \hline
        35 & 1 &  1   &440 439 \\ \hline
        43 & 1 &  1   &1 782 688 \\ \hline
        53 & 1 &  1   &505 767 \\ \hline
    \end{tabular}
    \end{center}
    \label{tab_agg1}
\end{table}

Now we apply the 2nd step of the aggregation procedure from Fig. \ref{fig_agg_steps} where we identify bases from the subsets, which are also indicated in Table \ref{tab_agg1}. As an example, from the 189 subsets of length 3, we identify that there are only 7 among them with different dual variables. Hence, there are 7 bases of length 3. We aggregate the 3-hour subsets that belong to the same basis, i.e., that have the exact same dual variables. Following the same procedure for each length in Table \ref{tab_agg1}, we found that there are only 70 bases in total. They have varying lengths ranging from 1 to 53 hours. In total, those 70 bases represent 681 hours out of the full 8736. When the bases are used in an aggregated model of reduced size\footnote{This represents a 92,2\% reduction in the number of time periods, and therefore variables, in the aggregated model}, they exactly approximate the complete model. This significantly reduces the computational burden\footnote{Especially if we take into account the exponential worst-case complexity of Simplex-like algorithms with respect to the number of variables and constraints \cite{Klee1972}.} of the model. 

\subsection{Results with Ramping and Network Constraints}
\label{ssec:resRampFlow}

We carry out the same data-aggregation procedure for the case with ramping and network constraints. The supporting tables are omitted due to space reasons, but can be found in \cite{BasisGithub}; most of the subsets lie in the 3 to 6-hour range with these being the ones that show the most significant aggregation in the 2nd Step. Moreover, the results also show that some extreme periods - even if they only occur once in the partition - have a significant impact on the objective function value. In summary, to completely represent this model, we only require 740 hours, 59 more than the previous 681.

A comparison of the three test cases: Network; Ramping; and, Network \& Ramping is presented in Table \ref{tab_summ_all}. \textit{Number of Bases} refers to how many different bases were needed to exactly approximate the full model. \textit{Max. Subset Length} indicates the maximum basis length (i.e., number of hours) in the aggregated model. \textit{Number Necessary Representative Hours} corresponds to how many hours, in total out of the 8736, are required by the Basis-Oriented approach to achieve zero error\footnote{Taking into account both, the number of bases and their length} in the aggregated model. Finally, \textit{Size Reduction: Variables} compares the size  of each aggregated model with the full one.

\subsection{Discussion}
\label{ssec:Discussion}

Our results highlight the complexity of temporal interlinking and show how the right-hand side value in one period, i.e., the ramping limit, affects the optimal solution, both forward and backward in time. Our research also shows that fixed-size representative periods might not be the best choice for aggregated models, and researchers must be ready to partition complete models in \textit{chunks} of variable temporal size, where each might represent a peculiarity of the system being analyzed; also the length of each subset does not necessarily correspond to the commonly used 24 or 168-hours periods. As our results show, wind production already challenges this convention, but a similar situation might entail if batteries, relying on efficient dynamic control algorithms, become widespread in the industrial and/or household sectors, as their consumption-charging cycles will not only depend on traditional variables like the time of the day but also on signaling from the environment (e.g., price). Finally, our results also support previous research about the importance of extreme values in an aggregated PSOM, as some of the periods that could not be merged happen only once or twice and have the highest objective-function values (e.g., subsets of length 13, 35 or 43 hours); hence the Basis-Oriented approach goes beyond the common assumptions of extreme days or weeks and tries to find these inseparable periods without imposing ex-ante conditions.

\begin{table}[htbp]
    \caption{Case Study Aggregation Comparison}
    \begin{center}
    \begin{tabular}{|c|c|c|c|c}
    \hline
                                & \textbf{Network} & \textbf{Ramping} & \makecell{\textbf{Network  \&}\\ \textbf{Ramping}} \\ \hline
    \textbf{\# of Bases}  & 5                & 70               & 85                          \\ \hline
    \textbf{Max. Subset Length} & 1                & 53               & 53                           \\ \hline
    \textbf{\# Required Repr. Hours} & 5                & 681              & 740                          \\ \hline
    \makecell{\textbf{Size Reduction:} \\ \textbf{Variables}}     & 99,94\%          & 92,20\%          & 91,53\%                      \\ \hline
    \end{tabular}
    \end{center}
    \label{tab_summ_all}
\end{table}

\section{Conclusion}
\label{sec:conclusion}

In this work, we addressed TSA when applied for optimization models, and in particular, for PSOMs. We extended Basis-Oriented TSA, which considers the structure of the optimization model in the aggregation process without losing accuracy, to include both network and ramping constraints. 

One of the contributions was to demonstrate that when the structure of temporal constraints is taken into account in the Basis-Oriented approach, an exact aggregation can still be obtained, and to develop a methodology to achieve this. To illustrate this, we developed a heuristic algorithm that performs an exploration of the dual space of the full PSOM solution. The algorithm identifies which hours should be grouped together, based on the constraints that are binding at any given time. In this manner, the length of the representative periods of the aggregated model is variable and arises naturally in an incremental way from the exploration procedure.

Another conclusion of our work is that if the goal is to perform temporal aggregation while keeping the results as close as possible to the complete solution, researchers must rethink their approach to using a uniform length in their representative periods (i.e., days) and start to consider aggregated models with different lengths that may or may not correspond to the intuitive time intervals commonly used in PSOM.

Finally, for future research, we plan to extend the Basis-Oriented approach to other types of constraints (i.e. including integrality such as unit commitment constraints), which were not addressed in this paper. Another topic of future research is inferring an aggregation directly from the input data and the optimization model's structure without already having a solution (primal and dual) of the full model.

\section*{Acknowledgment}

This work is part of the NetZero-Opt project, which has been funded via a Starting Grant of the European Research Council (ERC) (Grant agreement No. [101116212]).


%

\appendices
\section{Partitioning Algorithm}
\label{sec:appAlgorithm}

The algorithm is divided in two main parts composed by each of the loops; in the first one, the  ramping lengths of the model run are identified by associating each marginal cost with the closest integer multiple of the thermal unit's variable cost, $b$, the coefficient of $VC_{t}$ in \eqref{eq_marg_cost}. The second loop is the partitioning of the model's input data based on each period's marginal cost.

\begin{algorithm}[h]
\scriptsize
\caption{Dual-Based Partitioning}\label{alg:cap}
\begin{algorithmic}
\STATE $VC_t \gets $ \textit{Thermal generator variable cost}
\STATE $VC_w \gets $ \textit{Wind generator variable cost}
\STATE $VC_{nsp} \gets $ \textit{NSP cost}
\STATE $T \gets $ \textit{Ordered index set with the time representation}
\STATE $Basis[\forall t \in T] \gets \FALSE$
\STATE $Model \gets $ {\textit{Complete model solution}}
\STATE $MC[t] \gets $ \textit{Marginal costs from} $Model$
\STATE $DualRUP[t] \gets $ \textit{Ramp-up constraint duals from} $Model$
\STATE $DualRDN[t] \gets $ \textit{Ramp-down constraint duals from} $Model$
\STATE $DualWnd[t] \gets $ \textit{Wind availability duals from} $Model$

\FORALL{$mc \in MC | (mc > 0 \land mc \notin [VC_t, VC_w]) $}
    \STATE $Lengths[mc] \gets $ \textit{Multiple of} $VC_t$ closest to $ mc $
\ENDFOR

\FORALL{$t \in T$}
    \IF{$ (MC[t] \notin [VC_t, VC_w]) $}
        \IF{$MC[t] < 0$}
            \STATE Find next $t' > t$ such that $MC[t'] = VC_{nsp}$
            \STATE $Basis[t:t'] \gets \TRUE$
            \STATE $t \gets t' + 1$
        \ELSIF{$MC[t] = VC_{nsp}$}
            \STATE Find next $t' > t$ such that $MC[t'] \in [VC_t, VC_w]$
            \STATE $Basis[t:t'] \gets \TRUE$
        \ELSE
            \STATE $l \gets Lengths[MC[t]]$
            \IF{$DualRUP[t-1] > 0$}
                \STATE $t' \gets $ $max\{t'\} | (t' < (t-l) \land DualWnd[t'] > 0 )$
                \STATE $Basis[t':t] \gets \TRUE$
            \ELSIF{$DualRDN[t-1] > 0$}
                \STATE $t' \gets $ $min\{t'\} | (t' < (t+l) \land DualWnd[t'] > 0 )$
                \STATE $Basis[t:t'] \gets \TRUE$
            \ENDIF        
            \STATE $t \gets t' + 1$
        \ENDIF 
    \ENDIF

    \STATE $ Partition[p] \gets $ Subsets of longest contiguous $t$ \\ \quad \quad \quad such that $ Basis[t_{start}:t_{end}] = \TRUE$
\ENDFOR

\end{algorithmic}
\end{algorithm}

\ifCLASSOPTIONcaptionsoff
  \newpage
\fi



%



\bibliographystyle{IEEEtran}
\bibliography{bibtex/bib/IEEEabrv, bibtex/bib/references}
\vskip -2.5\baselineskip plus -1fil
%

\begin{IEEEbiographynophoto}{David Cardona-Vásquez}(Student Member, IEEE)
received a B.S. in computer engineering and M.S. in systems engineering degrees from the National University of Colombia in 2009 and 2011, respectively; and a B.S. in economics from EAFIT University, Colombia, in 2020. He is currently a PhD student at the Institute of Electricity Economics and Energy Innovation at Graz University of Technology, Austria.  His research interests include electricity markets, machine learning and optimization applications in economics, industrial organization and econometrics
\end{IEEEbiographynophoto}
\vskip -2\baselineskip plus -1fil

\begin{IEEEbiographynophoto}{Thomas Klatzer}(Student Member, IEEE)
obtained a Dipl.-Ing. (5-year degree) in Electrical Engineering and Business from Graz University of Technology. He is a PhD candidate at the Institute of Electricity Economics and Energy Innovation at Graz University of Technology, Austria. His research interests include renewable energy technologies, electricity markets and optimization of integrated power, gas, and hydrogen systems.
\end{IEEEbiographynophoto}

\vskip -2\baselineskip plus -1fil
\begin{IEEEbiographynophoto}{Sonja Wogrin}(Senior Member, IEEE)
holds a Dipl. Ing. (5-year degree) in Technical Mathematics from Graz University of Technology, a Master of Science from MIT, and a PhD in Power Systems from Comillas Pontifical University. Currently, she is the Head of the Institute of Electricity Economics and Energy Innovation at Graz University of Technology. Her research interests are centered around the field of modeling and optimization within the energy sector.

\end{IEEEbiographynophoto}



\end{document}